\documentclass[11pt]{article}
\usepackage{my_preamble}
\usepackage[small]{titlesec}
\usepackage[a4paper, margin=1in]{geometry}

\title{\textbf{The Hitchin morphism for certain surfaces fibered over a curve}}
\author{Matthew Huynh}
\date{\today}

\begin{document}
\maketitle

\begin{abstract}
    The Chen-Ng\^o Conjecture \cite{chen_ngo20}*{Conjecture 5.2} predicts that the Hitchin morphism from the moduli stack of $G$-Higgs bundles on a smooth projective variety surjects onto the space of spectral data (recalled in \S \ref{sect:hitchin_bases}).
    The conjecture is known to hold for the group $\rm{GL}_n$ and any surface \cite{song_sun_hitchin_morphism}, and for the group $\rm{GL}_2$ and any smooth projective variety \cite{he_liu_hitchin_morphism}.
    We prove the Chen-Ng\^o Conjecture for any reductive group when the variety is a ruled surface or (a blowup of) a nonisotrivial elliptic fibration with reduced fibers.
    Furthermore, if the group is a classical group, i.e. $G \in \{\rm{SL}_n,\rm{SO}_n,\rm{Sp}_{2n}\}$, then we prove that the Hitchin morphism restricted to the Dolbeault moduli space of semiharmonic $G$-Higgs bundles surjects onto the space of spectral data.

\end{abstract}

\tableofcontents

\section{Introduction}\label{sect:intro}
Let $G$ be a reductive group. The moduli stack of $G$-Higgs bundles on a smooth projective curve has a rich geometry and topology. 
N. Hitchin defined a surjective morphism from this moduli space to affine space \cite{hitchin_integrable_systems}, and the study of this morphism has produced stunning applications like Ng\^o's proof of the Fundamental Lemma in the Langlands program \cite{ngo_fund_lemma}.

We can replace the curve by a higher-dimensional smooth projective variety and obtain a moduli stack of $G$-Higgs bundles equipped with a morphism to an affine space.
This morphism, called the Hitchin morphism, is not necessarily surjective once the dimension of our variety is at least two.
However, there is a conjectured image $\mc{B}(X,G)$ defined by T. H. Chen and B. C. Ng\^o \cite{chen_ngo20}*{Conjecture 5.2}; its definition is recalled in \S \ref{sect:hitchin_bases}. 
Roughly speaking, the conjectured image $\mc{B}(X,G)$ is the locus where one can construct, for each $k$-point,  a ``cameral cover'' of the original variety and hope to obtain a modular description of the fiber of the Hitchin morphism in terms of the cameral cover.

If $X$ is an abelian variety of dimension $d$, then the Hitchin morphism for $X$ surjects onto $\mc{B}(X,G)$, essentially because the sheaf $\Omega^1_X$ is free of rank $d$ \cite{chen_ngo20}*{Example 5.1}. 
All other known results are for the general linear group, because in this case one may work with ``spectral covers'', which are simpler than cameral covers. 
When $X$ is a surface, Chen and Ng\^o prove a ``spectral correspondence'' result over an open subset $\mc{B}^{\heartsuit} \subset \mc{B}(X,\rm{GL}_n)$. 
More precisely, for every $b \in \mc{B}^{\heartsuit}(k)$, they construct a finite, flat cover $X_b^{\rm{CM}} \to X$, and prove that the fiber of the Hitchin morphism over $b$ is isomorphic to the stack of Cohen-Macaulay sheaves of generic rank 1 on $X_b^{\rm{CM}}$ \cite{chen_ngo20}*{Theorem 7.3}.
In particular, this shows that $\mc{B}^{\heartsuit}$ is contained in the image of the Hitchin morphism. 

However, the open subset $\mc{B}^{\heartsuit}$ is mysterious in general, and $\mc{B}^{\heartsuit}$ may not be dense in $\mc{B}(X,\rm{GL}_n)$, as the scheme $\mc{B}(X,\rm{GL}_n)$ need not be irreducible. 
L. Song and H. Sun prove that when $X$ is a surface, the conjectured image $\mc{B}(X,\rm{GL}_n)$ \textit{is} the image of the Hitchin morphism by ``decomposing'' spectral data $b \in \mc{B}(X,\rm{GL}_n)\setminus\mc{B}^{\heartsuit}$ and using the spectral correspondence \cite{song_sun_hitchin_morphism}.
S. He and J. Liu prove that for $G = \rm{GL}_2$ and any smooth projective variety $X$, the conjectured image $\mc{B}(X,\rm{GL}_2)$ \textit{is} the image of the Hitchin morphism. 
Roughly speaking, they construct Cohen-Macaulay covers $\tilde{X}_b$ of $X$ for all $b$ in an open subset of $\mc{B}(X,\rm{GL}_2)$, prove a spectral correspondence for these covers, and produce Higgs bundles over the remainder of $\mc{B}(X,\rm{GL}_2)$ explicitly \cite{he_liu_hitchin_morphism}.

For $G= \rm{GL}_n$ and $X$  a surface fibered over a curve $C$, Chen and Ng\^o show that the Hitchin base for the curve embeds as a closed subscheme in $\mc{B}(X,\rm{GL}_n)$, and that the image of the Hitchin morphism contains this subscheme \cite{chen_ngo20}*{\S 8}.

The purpose of this paper is to study the Hitchin morphism for surfaces fibered over a curve and \textit{other} reductive groups. 
The results are stated precisely in \S \ref{sect:summary_of_results} after establishing some notation in \S \ref{sect:notation}.
We prove that the Hitchin base for the curve embeds as a closed subscheme in the conjectured image $\mc{B}(X,G)$, and that the image of the Hitchin morphism contains this closed subscheme. 
When $G$ is a classical group and the fibered surface has only reduced fibers, we can use the standard linear representation of $G$ to deduce more, namely that the locus of semiharmonic $G$-Higgs bundles surjects onto the Hitchin base for the curve, which is embedded in $\mc{B}(X,G)$.
Finally, for ruled surfaces and blowups of nonisotrivial elliptic fibrations with reduced fibers, we prove the Chen-Ng\^o Conjecture as a corollary of these results. 

\subsection{Acknowledgments}
The author would like to thank Mark Andrea de Cataldo for helpful mathematical conversations and for useful comments on earlier drafts of this paper. 
He would also like to thank the anonymous referee for many helpful comments.
The author was supported by NSF Grant DMS-2200492.

\subsection{Notation and conventions}\label{sect:notation}
\begin{itemize}
    \item Let $k$ be an algebraically closed field of characteristic zero. 
    All geometric objects in this paper live over $\Spec(k)$. 
    \item In \S \ref{sect:higgs_stuff} and \S \ref{sect:pullback_higgs_bdles}, we use the Einstein summation convention when defining various maps. 
    \item Suppose that $G$ is an affine algebraic group, $Y$ is a scheme with a (left) $G$-action, $X$ is a scheme, and $P$ is a (right) principal $G$-bundle over $X$. 
    Then we denote the associated bundle $[P\times Y/G]$ by the symbol $P(Y)$. 
    If $E$ is a geometric vector bundle on $Y$, and $\rm{Fr}(E)$ is the associated principal $\rm{GL}_n$-bundle of frames, then we set
    \[
        E(Y) := [\rm{Fr}(E) \times Y/\rm{GL}_n].
    \]
    \item Let $X$ be a smooth, connected, projective variety of dimension $d$, and let $G$ be a reductive group. 
    Then $\mc{M}(X,G)$ denotes the moduli stack of $G$-Higgs bundles on $X$ (\S \ref{sect:higgs_moduli_space}), $\mc{A}(X,G)$ denotes the Hitchin base (\S \ref{sect:hitchin_bases}), $h_{X,G}$ denotes the Hitchin morphism $\mc{M}(X,G) \to \mc{A}(X,G)$ (\S \ref{sect:hitchin_morphism}), $\mc{B}(X,G) \subset \mc{A}(X,G)$ denotes the conjectured image of $h_{X,G}$ (\S \ref{sect:hitchin_bases}), and $M_{\rm{Dol}}(X,G)$ denotes the moduli space of semiharmonic $G$-Higgs bundles (\S \ref{sect:dolbeault_moduli_space}).
    \item A fibered surface $f:X \to C$ is shorthand for 
    \begin{enumerate}[label=(FS)]
        \item \label{hyp:fs} a smooth, connected, projective surface $X$, a smooth, connected, projective curve $C$, and a morphism $f: X \to C$ that is proper, flat, surjective and whose generic fiber is a projective, smooth, connected curve.
    \end{enumerate}
\end{itemize}

\subsection{Summary of contents and results}\label{sect:summary_of_results}
We recall definitions of objects appearing in the theorem statements and relevant preliminaries in section \ref{sect:preliminaries}. Nothing contained therein is new.
Section \ref{sect:pullback_higgs_bdles} contains all details of the proof of Theorem \ref{thm:red_gp_G}, and section \ref{sect:hitch_morphism_fibered_sf_classical_gp} contains all details of the proof of Theorem \ref{thm:classical_gp}.

Let $f: X \to C$ be a fibered surface \ref{hyp:fs}. Then for each positive integer $i$, the pullback of symmetric differentials $H^0(C,S^i\Omega^1_C) \to H^0(X,S^i\Omega^1_X)$ is injective, and thus defines a closed embedding $\mc{A}(C,G) \to \mc{A}(X,G)$. 
\begin{thm}
    \label{thm:red_gp_G}
    Let $f: X \to C$ be a fibered surface \ref{hyp:fs}, and let $G$ be a reductive group. 
    Then
    \begin{enumerate}
        \item the closed immersion $\mc{A}(C,G) \to \mc{A}(X,G)$ induced by pullback of symmetric differentials factors through the conjectured image of the Hitchin morphism $\mc{B}(X,G)$, and
        \item the image of the Hitchin morphism $h_{X,G}$ contains $\mc{A}(C,G)$. 
    \end{enumerate}
    \begin{proof}
        The proof of the first statement follows from Proposition \ref{prop:spaces_over_X}.
        
        A $G$-Higgs bundle $(E,\theta)$ on $C$ with $h_{C,G}(E,\theta) = a \in \mc{A}(C,G)$ pulls back to one on $X$, and this pullback $G$-Higgs bundle maps to $a$ via $h_{X,G}$ (Corollary \ref{cor:pullback_higgs_bdle_spec_data}).
        The second assertion in the theorem follows from the fact that the Hitchin morphism for $C$ is surjective. 
    \end{proof}
\end{thm}

\begin{thm}
    \label{thm:classical_gp}
    Let $f: X \to C$ be a fibered surface \ref{hyp:fs} with only reduced fibers, and let $G$ be a classical group, i.e. $G \in \{\rm{SL}_n,\rm{SO}_n,\rm{Sp}_{2n}\}$.
    Then the image of the Hitchin morphism restricted to the Dolbeault moduli space $M_{\rm{Dol}}(X,G)$ contains the Hitchin base $\mc{A}(C,G)$ of the curve. 
    \begin{proof}
        The Hitchin base $\mc{A}(C,G)$ for the curve embeds as a closed subscheme of $\mc{B}(X,G)$ by Theorem \ref{thm:red_gp_G}.

        If $C = \P^1$, then $\mc{A}(C,G)$ is a point, which maps to $0_{\mc{B}} \in \mc{B}(X,G)$. 
        The trivial $G$-bundle on $X$ with zero Higgs field satisfies the following properties: (i) $h(G\times X,0) = 0_{\mc{B}}$, (ii) the associated vector bundle is trivial, and hence has zero Chern classes, and (iii) the associated vector bundle is slope-semistable, as it is the direct sum of slope-stable sheaves that are all of the same slope.
        Therefore, this $G$-Higgs bundle is semiharmonic.

        Now assume that $g(C) \geq 1$.
        The restriction of the Hitchin morphism to the Dolbeault moduli space $M_{\rm{Dol}}(X,G)$ is proper, so it suffices to produce a nonempty open subset of $\mc{A}(C,G)$ contained in the image of $h_{X,G}$ restricted to $M_{\rm{Dol}}(X,G)$. (The Hitchin base $\mc{A}(C,G)$ is irreducible, so any nonempty open set is dense). 

        For all $a \in \mc{A}(C,G)$, there is a companion $G$-Higgs bundle $(E_a',\theta_a')$ on $C$ with $h_{C,G}(E'_a,\theta'_a) = a$, (see \cite{hameister_ngo} for instance).
        For generic $a \in \mc{A}(C,G)$, this companion $G$-Higgs bundle pulls back to a semiharmonic $G$-Higgs bundle on $X$, i.e. a $k$-point of $M_{\rm{Dol}}(X,G)$ (Proposition \ref{prop:change_of_group}, \S \ref{sect:hitch_morphism_fibered_sf_classical_gp}).
        The theorem follows.
    \end{proof}
\end{thm}

When $f: X \to C$ is a ruled surface or nonisotrivial elliptic fibration with only reduced fibers, then the pullback of symmetric differentials $H^0(C,S^i\Omega^1_C) \to H^0(X,S^i\Omega^1_X)$ is an isomorphism for all $i$ \cite{chen_ngo20}*{Proposition 8.1}.
This implies that $\mc{A}(C,G) = \mc{B}(X,G) = \mc{A}(X,G)$, so by Theorems \ref{thm:red_gp_G} and \ref{thm:classical_gp} we deduce the following 
\begin{cor}
    \label{cor:conj_for_special_fibered_sfs}
    Let $f: X \to C$ be a ruled surface or nonisotrivial elliptic fibration with only reduced fibers, and let $G$ be a reductive group. 
    Then the Hitchin morphism surjects onto $\mc{A}(X,G) = \mc{B}(X,G)$.
    If $G$ is a classical group, then the Hitchin morphism restricted to the Dolbeault moduli space $M_{\rm{Dol}}(X,G)$ surjects onto $\mc{A}(X,G) = \mc{B}(X,G)$. 
\end{cor}

Suppose that $X$ is a smooth projective surface, and $\pi: Y \to X$ is a birational morphism. 
Then the pullback of symmetric differentials $H^0(X,S^i\Omega^1_X) \to H^0(Y,S^i\Omega^1_Y)$ is an isomorphism for all $i$ (\S \ref{sect:pullback_symm_diff_sfs}).
Thus we can determine the image of the Hitchin morphism for a slightly larger class of surfaces:
\begin{cor}
    Let $X$ be a nonisotrivial elliptic fibration with only reduced fibers, let $\pi: Y \to X$ be a birational morphism, and let $G$ be a reductive group. 
    Then the Hitchin morphism $h_{Y,G}: \mc{M}(Y,G) \to \mc{A}(Y,G)$ surjects onto $\mc{A}(Y,G) = \mc{B}(Y,G)$.
    \begin{proof}
        Without loss of generality, suppose that $\pi: Y \to X$ is the blowup of a point.
        The pullback of symmetric differentials $H^0(X,S^i\Omega^1_X) \to H^0(Y,S^i\Omega^1_Y)$ identifies the Hitchin bases for $X$ and $Y$, i.e. we have $\mc{A}(X,G) = \mc{A}(Y,G)$ (\S \ref{sect:pullback_symm_diff_sfs}).
        The assumption on $X$ further implies that $\mc{A}(Y,G) = \mc{B}(Y,G)$.
        The pullback of a $G$-Higgs bundle $(E,\theta)$ on $X$ with $h_{X,G}(E,\theta) = a$ satisfies $h_{Y,G}(\pi^{\ast}E,\pi^{\ast}\theta) = a$ by Proposition \ref{prop:diagram_morphism_of_sfs}.
        We conclude by applying Corollary \ref{cor:conj_for_special_fibered_sfs}.
    \end{proof}
\end{cor}

\section{Preliminaries}\label{sect:preliminaries}
\subsection{Associated bundles}\label{sect:assoc_bdle}
Let $G$ be an affine algebraic group, $Y$ an affine scheme with a (left) $G$-action, $X$ a scheme, and $P \to X$ a (right) principal $G$-bundle. 
Then the associated bundle $P(Y)$ is affine over $X$; (the projection $P \times Y \to P$ is affine and $G$-equivariant, so the claim follows by descent). 

Thus the pushforward of $\mc{O}_{P(Y)}$ to $X$, denoted $\mc{A}$, is an $\mc{O}_X$-algebra. 
Suppose further that $P$ is Zariski-locally trivial, and that $\{U_{\alpha} \to X\}_{\alpha \in I}$ is an affine open cover of $X$ over which $P$ trivializes. 
Denote the transition functions by $g_{\alpha\beta} \in G(\Gamma(U_{\alpha\beta},\mc{O}))$ and say $Y = \Spec(A)$. 
Then $\mc{A}|_{U_{\alpha}} = (R_{\alpha}\otimes_k A)^{\sim}$, and these sheaves are glued together on affine overlaps $V = \Spec(R) \subset U_{\alpha\beta}$ by the $R$-module isomorphisms
\begin{equation*}
    \label{eq:trans_map_assoc_bdle}
    \xymatrixcolsep{1.5cm}\xymatrix{
        R \otimes_k A \ar[r]^-{\id \otimes\, \rm{coact}} & R \otimes_k A \otimes_k k[G] \ar[r]^-{\rm{id} \otimes\, g_{\alpha\beta}} & R \otimes_k A \otimes_k R \ar[r]^-{m_{13}} & R \otimes_k A.
    }
\end{equation*}

Observe that the associated bundle construction behaves well with respect to pullback. 
In other words, if $f: X' \to X$ is a morphism of $k$-schemes, then we have an isomorphism of $X'$-schemes $(f^{\ast}P)(Y) \cong P(Y) \times_X X'$. 
To verify this claim, observe that there exists a commutative diagram,
\begin{equation*}
    \label{eq:assoc_bdle_pullback}
    \xymatrix{
        f^{\ast}P \times Y \ar[r] \ar[d] & f^{\ast}P \ar[d] \ar[r] & X' \ar[d]^{f} \\
        P \times Y \ar[r] & P \ar[r] & X,
    }
\end{equation*}
where the horizontal maps are projections, and both squares are cartesian. 
The arrows in the left square are all $G$-equivariant, so the claim follows by descent.

\subsection{Higgs bundles and the Hitchin morphism}\label{sect:higgs_stuff}
In this subsection, we review the definitions of the moduli stack of $G$-Higgs bundles, the Hitchin base, the Hitchin morphism, and its conjectured image in terms of mapping stacks. 
Everything presented from \S \ref{sect:comm_scheme} to \ref{sect:hitchin_morphism} is from \cite{chen_ngo20}*{\S 3-5}.

\subsubsection{The commuting scheme}\label{sect:comm_scheme}
Let $G$ be a reductive group, and let $\mf{g} = \rm{Lie}(G)$. 
Fix a maximal torus $T \subset G$, let $\mf{t} = \rm{Lie}(T)$, and let $r$ be the rank of $G$, (i.e. $\dim T$). 
Let $W = N_G(T)/T$ be the Weyl group. 
Choose homogeneous generators $c_1,\ldots,c_n$ of $k[\mf{g}]^G \cong k[\mf{t}]^W$, and denote their degrees by $e_1,\ldots,e_n$ respectively.

For any positive integer $d$, we can define the $d$-fold commutator map
\begin{equation*}
    \label{eq:commutator_map}
    \mf{g}^d \to \prod_{1\leq i < j \leq d} \mf{g},\quad\quad (\theta^1,\ldots,\theta^d) \mapsto ([\theta^i,\theta^j])_{i < j}.
\end{equation*}
We define the $d$-fold commuting scheme $\mf{C}^d$ as the scheme-theoretic fiber over $(0,\ldots,0) \in \prod_{i < j}\mf{g}$. Note that $\mf{C}^1 = \mf{g}$.

The group $G$ acts on $\mf{g}^d$ by the diagonal adjoint action, and this action restricts to one on $\mf{C}^d$. 
There is also an action of $\rm{GL}_d$ on $\mf{g}^d$ defined by
\begin{equation*}
    \label{eq:GL_d_action_on_comm_scheme}
    \rm{GL}_d(R) \times \Hom_k((\mf{g}^d)^{\ast},R) \to \Hom_k((\mf{g}^d)^{\ast},R), \quad\quad
    (x^i_j),(\theta^1,\ldots,\theta^d) \mapsto (x^1_j\theta^j,\ldots,x^d_j\theta^j),
\end{equation*}
where $R$ is any $k$-algebra. 
This action on $\mf{g}^d$ also restricts to one on $\mf{C}^d$ because the Lie bracket is bilinear, and the $G$ and $\rm{GL}_d$ actions commute. When $d = 1$, the $\rm{GL}_1 = \G_m$ action on $\mf{C}^1 = \mf{g}$ is just scaling.

\subsubsection{The moduli stack of Higgs bundles}\label{sect:higgs_moduli_space}
Let $X$ be a connected, smooth, projective variety of dimension $d$, and let $T^{\ast}_X = \uSpec(S^{\bullet}(\Omega^1_X)^{\vee})$ be its cotangent bundle. 
By the preceding discussion, we can form the associated bundle $T^{\ast}_X(\mf{C}^d)$. 
Notice that $G$ acts on $T^{\ast}_X(\mf{C}^d)$, (the action is the fiberwise diagonal adjoint action of $G$ on $\mf{C}^d$), and the morphism $T^{\ast}_X(\mf{C}^d) \to X$ is $G$-invariant. 

The moduli stack of $G$-Higgs bundles on $X$ is defined as
\begin{equation*}
    \label{eq:higgs_moduli_stack}
    \mc{M}(X,G) := \rm{Maps}_X(X,[T^{\ast}_X(\mf{C}^d)/G]).
\end{equation*}
A $G$-Higgs bundle on $X$, (i.e. $k$-point of this stack), is the data of a principal $G$-bundle $E \to X$, and an $\mc{O}_X$-module morphism $\theta: (\Omega^1_X)^{\vee} \to \rm{ad}(E)$ such that $[\theta(u),\theta(v)] = 0$ for all local sections $u,v$ of $(\Omega^1_X)^{\vee}$. 
This last condition is called the \textit{integrability condition}, and is typically abbreviated as $\theta\wedge\theta = 0$. 

A $\rm{GL}_n$-Higgs bundle on $X$ is equivalent to the data of a locally free sheaf of rank $n$, say $E$, and an $\mc{O}_X$-module morphism $\theta: (\Omega^1_X)^{\vee} \to \rm{End}(E)$ satisfying the integrability condition. 
We use this identification implicitly throughout the paper. 

\subsubsection{The Hitchin base and conjectured image of the Hitchin morphism}\label{sect:hitchin_bases}
Here we deviate slightly from the more intrinsic descriptions of the Hitchin base and conjectured image given in \cite{chen_ngo20}*{\S 4}. 
We present the constructions following \cite{losik2006polarizations} to facilitate later computations. 

Let $V$ denote the $k$-vector space $k^d$, and let $v_1,\ldots,v_d$ denote the standard basis. 
Let $\mf{A}$ be the scheme associated to the product of the vector spaces $S^{e_j}V$ for $j = 1,\ldots,n$; when $d = 1$, we have $\mf{A} = \mf{c} := \Spec(k[\mf{g}]^G)$. 
Recall that a $d$-tuple of non-negative integers $\ul{i} := (i_1,\ldots,i_d)$ whose sum equals $e$ is called a \textit{weak composition of} $e$ \textit{of length} $d$. 
The set of vectors 
\[
    \{v^{j,\ul{i}} := v_1^{i_1}\cdots v_d^{i_d} \big|\, \text{$\ul{i}$ is a weak composition of $e_j$ of length $d$}\,\}, 
\]
forms a basis of $S^{e_j}V$. 
Denote by $z_{j,\ul{i}}$ the dual basis. 
Then the ring of regular functions on $\mf{A}$ is a polynomial ring generated by the symbols $z_{j,\ul{i}}$. 

Note that the standard representation of $\rm{GL}_d$ on $V$ induces a $\rm{GL}_d$-action on $\mf{A}$. 
Let $X$ be as in \S \ref{sect:higgs_moduli_space}. 
The Hitchin base (for $X$) is defined as a scheme of sections,
\begin{equation*}
    \label{eq:hitchin_base}
    \mc{A}(X,G) := \rm{Maps}_X(X,T^{\ast}_X(\mf{A})).
\end{equation*}
The $X$-scheme $T^{\ast}_X(\mf{A})$ is isomorphic to the geometric vector bundle associated to the locally free sheaf $\oplus_{j = 1}^n S^{e_j}\Omega^1_X$, so $\mc{A}(X,G)$ is isomorphic to affine space. 

Next, we describe the conjectured image of the Hitchin morphism, also called the \textit{space of spectral data}. There is a subring of $k[\mf{t}^d]^W$ called the \textit{polarization ring}. 
To define it, consider the map of vector spaces
\begin{equation*}
    \label{eq:polarization_map}
    \varphi: \mf{t}^d \times V \to \mf{t},\quad\quad (\theta^1,\ldots,\theta^d),(b_1,\ldots,b_d) \mapsto b_i\theta^i.
\end{equation*}
For any $W$-invariant function $f \in k[\mf{t}]^W$, there is a decomposition of $\varphi^{\ast}f$ into isotypic components for the scaling action of $\G_m$ on $V$:
\begin{equation}
    \label{eq:polarization_decomp}
    \varphi^{\ast}(f)(\theta^1,\ldots,\theta^d;b_1,\ldots,b_d) = \sum_{(i_1,\ldots,i_d) \in \Z^d_{\geq 0}} b_1^{i_1}\cdots b_d^{i_d} f_{(i_1,\ldots,i_d)}(\theta^1,\ldots,\theta^d).
\end{equation}
The functions $f_{(i_1,\ldots,i_d)}$ are called the polarizations of $f$. 
The polarization subring of $k[\mf{t}^d]^W$, denoted $k[\mf{t}^d]^{\rm{pol}}$, is by definition the ring generated by all polarizations of $W$-invariant functions on $\mf{t}$. 

If $f$ is a homogeneous polynomial function of degree $e$, then the only non-zero coefficients that appear in the decomposition (\ref{eq:polarization_decomp}) correspond to weak compositions of $e$ of length $d$. 
Furthermore, if $c_1,\ldots,c_n$ is a set of homogeneous polynomials that generate $k[\mf{t}]^W \cong k[\mf{g}]^G$, then the $\ul{i}$-polarizations of $c_j$, denoted $c_{j,\ul{i}}$, generate $k[\mf{t}^d]^{\rm{pol}}$. 
Let $\mf{B}$ be the spectrum of the ring of polarizations. 

\begin{ex}
    We describe the scheme $\mf{B}$ in the case when $G = \rm{SL}_2$ and $d = 2$. 
    The description here complements the intrinsic description of this example presented in \cite{chen_ngo20}*{Example 4.2}.

    The quotient $\mf{t} \gitq W$ is isomorphic to the affine line, and the map $\mf{t} \to \mf{t} \gitq W$ is given by sending the coordinate $c$ on $\mf{t} \gitq W$ to the polynomial function $\det$ on $\mf{t}$. 
    Next we compute the polarizations of $\det$. 
    Given $\theta^i = \rm{diag}(a_i,-a_i) \in \mf{t}(k)$ and $b_i \in \A^1(k)$, we see that
    \[
        \det(b_i\theta^i) = -(b_1^2a_1^2 + 2b_1b_2a_1a_2 + b_2^2a_2^2),
    \]
    so the polarizations of the function $c = \det$ are
    \[
        c_{(2,0)} = -\det(\theta^1),\quad c_{(1,1)} = -\rm{tr}(\theta^1\theta^2),\quad c_{(0,2)} = -\det(\theta^2).
    \]
    The scheme $\mf{B}$ is equal to the spectrum of the ring
    \[
        k[c_{(2,0)},c_{(1,1)},c_{(0,2)}]/(c_{(1,1)}^2 - 4c_{(2,0)}c_{(0,2)}).
    \]

    We also have $\mf{t}^2\gitq W \simeq \mf{B}$, as we now explain. 
    The nontrivial element of $W = \mf{S}_2$ acts on $\mf{t}^2 := \Spec(k[x,y])$ by $x \mapsto -x$ and $y \mapsto -y$. 
    The polynomials $x^2,2xy,y^2$ generate the ring of $W$-invariant functions on $\mf{t}^2$, and it is clear that these generators satisfy the same relation satisfied by the polarizations of $c = \det$. One also sees in this example that the morphism $\mf{t}^d \to \mf{t}^d \gitq W$ is not flat in general when $d \geq 2$. 
\end{ex}

The $\rm{GL}_d$-action on $\mf{C}^d$ restricts to $\mf{t}^d \subset \mf{C}^d$, and this induces a $\rm{GL}_d$-action on $\mf{B}$.
The conjectured image is defined as a scheme of sections,
\begin{equation*}
    \label{eq:space_of_spec_data}
    \mc{B}(X,G) := \rm{Maps}_X(X,T^{\ast}_X(\mf{B})).
\end{equation*}

There is a closed immersion $\mf{B} \to \mf{A}$ defined by the ring map $z_{j,\ul{i}} \to c_{j,\ul{i}}$ that respects the $\rm{GL}_d$-actions on both schemes. 
This induces a closed immersion $\mc{B}(X,G) \subset \mc{A}(X,G)$.

\subsubsection{The Hitchin morphism for curves and surfaces}\label{sect:hitchin_morphism}
Although the Hitchin morphism is defined for smooth projective varieties of all dimensions, for neatness of exposition, we only describe it for dimensions 1 and 2. 
The descriptions below also agree with the definition in \cite{simpson_moduli2} when $G$ is the general linear group.

Suppose $C$ is a smooth projective curve. 
The adjoint quotient map $\chi: \mf{g} \to \mf{c}$ induces a morphism $[\chi]: [T^{\ast}_C(\mf{g})/G] \to T^{\ast}_C(\mf{c})$, and the Hitchin morphism $h_{C,G}: \mc{M}(C,G) \to \mc{A}(C,G)$ is the induced morphism of mapping stacks. 
If we choose a non-zero root vector $f \in \mf{g}_{-\alpha_i}$ for each simple root $\alpha_i$, and a square root of the canonical bundle of $C$, then we obtain a section of the Hitchin fibration, called the Kostant section (more details on this construction can be found in \cite{Chen_Zhu_2017}*{\S 2.3}). 

Suppose $X$ is a smooth projective surface. 
Then the restriction map $k[\mf{C}^2]^G \to k[\mf{t}^2]^W$ is an isomorphism \cite{li_nadler_yun}. 
Thus we obtain a composition $\mf{C}^2 \to \mf{B} \to \mf{A}$ that is $\rm{GL}_2$-equivariant and $G$-invariant, which induces a composition of $X$-morphisms 
\begin{equation*}
    [T^{\ast}_X(\mf{C}^2)/G] \to T^{\ast}_X(\mf{B}) \to T^{\ast}_X(\mf{A}).
\end{equation*}
The Hitchin morphism is the induced morphism of mapping stacks $h_{X,G}: \mc{M}(X,G) \to \mc{A}(X,G)$. 

By construction, the Hitchin morphism $h_{X,G}$ always factors through $\mc{B}(X,G)$.
Although we have not described $h_{X,G}$ when $\dim(X) \geq 3$, it is still true that the image of any geometric point $(E,\theta) \in \mc{M}(X,G)(k)$ lies in $\mc{B}(X,G)(k)$ \cite{chen_ngo20}*{Proposition 5.1}.
If $(E,\theta)$ is a $G$-Higgs bundle on $X$ that maps to $b \in \mc{B}(X,G)(k)$, then we call $b$ the spectral data of $(E,\theta)$. 

\subsubsection{Dolbeault moduli spaces}\label{sect:dolbeault_moduli_space}
Introducing a notion of semistability allows us to obtain good moduli spaces in the sense of \cite{alper_gms}.
In particular, there are schemes parametrizing semistable Higgs bundles with vanishing Chern classes constructed by Simpson in \cites{simpson_moduli1, simpson_moduli2}. 
We recall the salient details now. 

Let $X$ be a connected, smooth, projective variety of dimension $d$, and fix an ample line bundle $L$ on $X$. 
A Higgs sheaf $(E,\theta)$ on $X$ is the data of a coherent sheaf $E$ and a morphism of $\mc{O}_X$-modules $\theta: E \to E \otimes \Omega^1_X$ satisfying $\theta\wedge\theta = 0$, i.e. the associated morphism $(\Omega^1_X)^{\vee} \to \rm{End}(E)$ satisfies the integrability condition described in \S \ref{sect:higgs_moduli_space}. 
A sub-Higgs sheaf of $E$ is a subsheaf $F \subset E$ such that $\theta|_F$ factors through $F \otimes \Omega^1_X \to E \otimes \Omega^1_X$. 

We call a Higgs sheaf $(E,\theta)$ $p$-semistable if $E$ is a pure $d$-dimensional coherent sheaf, and for all sub-Higgs sheaves $F \subset E$ with $0 < \rm{rk}(F) < \rm{rk}(E)$, there exists an integer $N > 0$ such that the inequality
\begin{equation*}
    \frac{P(F,n)}{\rm{rk}(F)} \leq \frac{P(E,n)}{\rm{rk}(E)}
\end{equation*}
holds for all $n \gg N$, where $P(E,n)$ is the Hilbert polynomial of $E$ with respect to $L$. 
Similarly, we call a Higgs sheaf $\mu$-semistable it is pure of dimension $d$, and the inequality $\mu(F) \leq \mu(E)$ holds for all sub-Higgs sheaves of strictly smaller rank, (recall $\mu = \rm{deg}/\rm{rk}$, and $\rm{deg}(E) = c_1(E)\cdot L^{d-1}$ because $X$ is smooth and projective). 
Both semistability conditions can be checked via quotients, that is, $(E,\theta)$ is $p$-semistable (resp. $\mu$-semistable) if and only if for all proper, purely $d$-dimensional quotients $q: E \to F'$ with $\ker(q)$ a sub-Higgs sheaf, we have $P(E,n)/\rm{rk}(E) \leq P(F',n)/\rm{rk}(F')$, (resp. $\mu(E) \leq \mu(F')$).

Now let $G$ be a reductive group. 
We say that a $G$-Higgs bundle $(E,\theta)$ is semiharmonic (or is of semiharmonic type), if the Chern classes of $E$ vanish in rational cohomology, and if there exists a faithful representation $\rho: G \to \rm{GL}(V)$ such that the associated Higgs bundle is $p$-semistable. 

Fix a point $x \in X$, and consider the following moduli functor
\begin{equation*}
    \label{eq:representation_space}
    \begin{gathered}
        R^{\natural}_{\rm{Dol}}(X,x,G): (\rm{Sch}/k)^{\rm{op}} \to \rm{Sets} \\
        S \mapsto \left\{\parbox[h]{9cm}{triples \ensuremath{(E,\theta,b)} where \ensuremath{(E,\theta)} is a \ensuremath{G}-Higgs bundle on \ensuremath{X_S}, such that for each closed point \ensuremath{s \in S}, \ensuremath{(E_s,\theta_s)} is of semiharmonic type, and \ensuremath{b: S \to E|_{\{x\}\times S}} is a section}\right\}/\text{isomorphism}.
    \end{gathered}
\end{equation*}
Then there is a scheme $R_{\rm{Dol}}(X,x,G)$ that represents this moduli functor \cite{simpson_moduli1}*{Theorem 4.10}. 
Furthermore, if $G \to H$ is a closed embedding, then there is a closed embedding $R_{\rm{Dol}}(X,x,G) \to R_{\rm{Dol}}(X,x,H)$. 

The group $G$ acts on the representation space $R_{\rm{Dol}}(X,x,G)$ by acting on the section. 
Simpson shows that one can form the GIT quotient for this action, and obtain a space $M_{\rm{Dol}}(X,G)$ that universally corepresents the functor sending a $k$-scheme $S$ to the set of isomorphism classes of $G$-Higgs bundles on $X_S$ that are semiharmonic when restricted to each fiber. 

We can instead take the stack quotient to obtain the moduli stack $\mc{M}_{\rm{Dol}}(X,G)$, which sends a $k$-scheme $S$ to the groupoid of $G$-Higgs bundles on $X_S$ that are semiharmonic when restricted to each fiber. 
There is a good moduli space morphism $\mc{M}_{\rm{Dol}}(X,G) \to M_{\rm{Dol}}(X,G)$ that comes from the quotient stack/GIT quotient construction \cite{alper_gms}*{Theorem 13.6, Remark 13.7}. 
Observe that $\mc{M}_{\rm{Dol}}(X,G)$ is a locally closed substack of $\mc{M}(X,G)$, because semiharmonicity is an open condition, and Chern classes in rational cohomology are locally constant in flat families.

The restriction of the Hitchin morphism to $\mc{M}_{\rm{Dol}}(X,G)$ factors through its good moduli space $M_{\rm{Dol}}(X,G)$, and the Hitchin morphism restricted to $M_{\rm{Dol}}(X,G)$ is proper (\cite{simpson_moduli2}*{Theorem 6.11} for the general linear group, and, for example, \cite{deCataldo_2021projective}*{Proposition 2.2.2} for reductive groups).
Therefore, if $U \subset \mc{B}(X,G)$ is any open subset contained in the image of $h_{X,G}$ restricted to $M_{\rm{Dol}}(X,G)$, then the closure of $U$ is also contained in the image. 

\subsection{Spectral covers}\label{sect:spectral_covers}
When $G$ is the general linear group $\rm{GL}_n$ and $C$ is a curve, there is a well-known description of the fibers of the Hitchin morphism in terms of spectral curves, which are finite, flat covers of $C$ \cites{BNR, schaub_98}. 
For now, we fix the group and the curve, so we omit these from the notation.

Let $a_i = (-1)^{i}\tr(\wedge^i -) \in k[\mf{gl}_n]^{\rm{GL}_n}$, let $A_n = k[a_1,\ldots,a_n]$, and let $\mf{c}_n = \Spec(A_n)$. 
Then define $B_n = A_n[x]/(x^n + a_1x^{n-1} + \cdots + a_n)$ and let $\mf{s}_n = \Spec(B_n)$. 
Observe that $B_n$ is a smooth $k$-algebra and a free $A_n$-module of rank $n$ equipped with an $A_n$-module endomorphism given by multiplication by $\ol{x}$. 
Therefore, we have a finite, flat morphism $p: \mf{s}_n \to \mf{c}_n$. 
Furthermore, we can extend the $\G_m$-action on $A_n$ to $B_n$, so that $p$ is $\G_m$-equivariant. 
Thus we can twist by $T^{\ast}_C$ to obtain a $C$-morphism $T^{\ast}_C(\mf{s}_n) \to T^{\ast}_C(\mf{c}_n)$. 

Given $a \in \mc{A}$, we obtain a finite, flat cover $p: C_a \to C$ by base change of $T^{\ast}_C(\mf{s}_n) \to T^{\ast}_C(\mf{c}_n)$, and $C_a$ is a closed subscheme of $T^{\ast}_C$. 
The spectral correspondence refers to the isomorphism between the fiber of the Hitchin morphism over $a$, and the moduli stack of torsion-free coherent sheaves of rank 1 on $C_a$. 
Torsion-free sheaves of rank 1 on $C_a$ pushforward to locally free sheaves of rank $n$ on $C$, and the Higgs field is given by the $\mc{O}_{T^{\ast}_C}$-module structure. 

\subsubsection{The companion section for the general linear group}\label{sect:companion_sect_gl}
In particular, given $a \in \mc{A}$ and a spectral cover $p: C_a \to C$, we can pushforward the structure sheaf to obtain the so-called \textit{companion Higgs bundle} $(E_a,\theta_a)$ on $C$. 
For every $a = (a_1,\ldots,a_n) \in \oplus_{i = 1}^n H^0(C,S^i\Omega^1_C)$, the locally free sheaf $E_a$ is isomorphic to $\oplus_{i = 0}^{n-1} (\Omega^1_C)^{-i}$, and the Higgs field $\theta_a$ is given by the companion matrix, 
\begin{equation*}
    \label{eq:companion_matrix}
    \mc{O}_C \oplus (\Omega^1_C)^{-1} \oplus \cdots \oplus (\Omega^1_C)^{-(n-1)} \to \Omega^1_C \oplus \mc{O}_C \oplus \cdots \oplus (\Omega^1_C)^{-(n-2)},\quad
        \begin{bmatrix}
        0 & 0 & \cdots & 0 & \otimes -a_n \\
        \id & 0 & \cdots & 0 & \otimes -a_{n-1} \\
        0 & \id & \cdots & 0 & \otimes -a_{n-2} \\
        \vdots & \vdots & \ddots & \vdots & \vdots \\
        0 & 0 & \cdots & \id & \otimes -a_1
    \end{bmatrix}.
\end{equation*}
This construction gives another section of the Hitchin morphism $h_{C,\rm{GL}_n}$, called the companion section. 

\subsubsection{The companion section for classical groups}\label{sect:companion_sect_classical_gp}
Let $G$ be a classical group and let $G \to \rm{GL}_N$ be the standard linear representation of $G$. 
Then $G$-Higgs bundles on a smooth projective variety $X$ can be described as $\rm{GL}_N$-Higgs bundles on $X$ with extra structure. (Recall we identify $\rm{GL}_N$-Higgs bundles with pairs $(E,\theta)$ where $E$ is a locally free sheaf of rank $N$, and $\theta$ is a Higgs field satisfying the integrability condition). 
\begin{enumerate}[label=($A_n$)]
    \item 
    An $\rm{SL}_{n+1}$-Higgs bundle is the data $(E,\theta,\tau)$, where $(E,\theta)$ is a $\rm{GL}_{n+1}$-Higgs bundle, and $\tau: \det(E) \isomto \mc{O}_X$ is a trivialization of the determinant line bundle of $E$, such that $\theta$ is ``traceless'', i.e. the projection of $h_{X,G}(E,\theta)$ to $H^0(X,\Omega^1_X)$ is zero. 
\end{enumerate}
\begin{enumerate}[label=($B_n$ and $D_n$)]
    \item 
    An $\rm{SO}_{2n+1}$-Higgs bundle (resp. $\rm{SO}_{2n}$-Higgs bundle) is the data $(E,\theta,\omega,\tau)$, where $(E,\theta)$ is a $\rm{GL}_N$-Higgs bundle, $\omega: S^2E \to \mc{O}_X$ is a fiberwise-nondegenerate symmetric pairing, and $\tau: \det(E) \isomto \mc{O}_X$ is a trivialization of the determinant line bundle of $E$, such that the Higgs field $\theta$ is anti-self-adjoint with respect to $\omega$, i.e. we must have $\omega(\theta(u),v) + \omega(u,\theta(v)) = 0$ in $\Omega^1_X$ for all local sections $u,v$ of $E$. 
\end{enumerate}
\begin{enumerate}[label=($C_n$)]
    \item 
    An $\rm{Sp}_{2n}$-Higgs bundle is the data $(E,\theta,\omega)$ where $(E,\theta)$ is a $\rm{GL}_N$-Higgs bundle, and $\omega: \Lambda^2 E \to \mc{O}_X$ is a fiberwise-nondegenerate symplectic pairing such that $\theta$ is anti-self-adjoint with respect to $\omega$. 
\end{enumerate}

For each classical group $G$, we may obtain a $G$-Higgs bundle using the companion $\rm{GL}_N$-Higgs bundle described in \S \ref{sect:companion_sect_gl}. 
The construction of these $G$-Higgs bundles dates back at least to Hitchin's seminal paper \cite{hitchin_integrable_systems}.
We call these \textit{companion $G$-Higgs bundles} and recount the construction below following the presentation of \cite{hameister_ngo}.

Let $A$ denote $k[\mf{g}]^G$ and let $A_N$ denote $k[\mf{gl}_N]^{\rm{GL}_N}$. 
Recall that $\mf{c}_N = \Spec(A_N)$ has a finite, flat, $\G_m$-equivariant cover $p: \mf{s}_N \to \mf{c}_N$, where $\mf{s}_N = \Spec(B_N)$ and $B_N = A_N[x]/(x^N + a_1x^{N-1} + \cdots + a_N)$.
We have a map $\mf{c} = \Spec(A) \to \mf{c}_N = \Spec(A_N)$, (which is a closed immersion when $G \neq \rm{SO}_{2n}$), and obtain a finite, flat, $\G_m$-equivariant cover $\mf{s} = \Spec(B) \to \mf{c}$ by base change of $p: \mf{s}_N \to \mf{c}_N$. 
Observe that $B = B_N \otimes_{A_N} A$ is a free $A$-module of rank $N$, with an $A$-module endomorphism $\theta$ given by multiplication by $\ol{x}$. 

\begin{enumerate}[label=($A_n$)]
    \item \label{item:spec_lin_gp} $G = \rm{SL}_{n+1}$.
    Recall that $A_{n+1} \to A$ is given by the quotient of $(a_1) \subset A_{n+1}$, so 
    \[ 
        B = B_{n+1} \otimes_{A_{n+1}} A = A[x]/(x^{n+1} + a_2x^{n-1} + \cdots + a_{n+1}). 
    \]
    Notice that $B$ is smooth over $k$ by the Jacobian criterion. 
    
    Let $a \in \mc{A}(C,\rm{SL}_{n+1}) \subset \mc{A}(C,\rm{GL}_{n+1})$ and consider the spectral curve $p: C_a \to C$. 
    The companion Higgs bundle $E_a = p_{\ast}\mc{O}_{C_a}$ satisfies $\det(E_a) = (\Omega^1_C)^{-n(n+1)/2}$. 
    If we choose a square root $\Omega'$ of $\Omega^1_C$, then $E_a' := E_a \otimes (\Omega')^{n}$ is a locally free sheaf of rank $n+1$ with trivial determinant. 
    After choosing some trivialization $\tau: \det(E_a') \isomto \mc{O}_C$, we obtain an $\rm{SL}_{n+1}$-Higgs bundle on $C$. 
\end{enumerate}
\begin{enumerate}[label=($B_n$)]
    \item \label{item:odd_spec_orth_gp} $G = \rm{SO}_{2n+1}$. 
    Recall that $A_{2n+1} \to A$ is given by the quotient of $(a_1,a_3,\ldots,a_{2n+1}) \subset A_{2n+1}$, so, 
    \[ 
        B = B_{2n+1} \otimes_{A_{2n+1}} A = A[x]/(x(x^{2n}+a_2x^{2n-2}+\cdots+a_{2n})).
    \]
    Notice that $B$ is not a smooth $k$-algebra, as one verifies by the Jacobian criterion. 

    There is a bilinear form $\omega: B \otimes_A B \to A$ that is nondegenerate, symmetric, and anti-self-adjoint with respect to $\theta$. 
    Furthermore, $\omega$ satisfies the following $\G_m$-compatibility:
    \begin{equation*}
        \label{eq:G_m_compat_odd_spec_orth_gp}
        (\omega \otimes \id)(\rm{coact}_{B\otimes_A B}(b_1\otimes b_2)) = \rm{coact}_A(\omega(b_1\otimes b_2))(1\otimes t^{2n}).
    \end{equation*}

    Let $a \in \mc{A}(C,G) \subset \mc{A}(C,\rm{GL}_{2n+1})$ and consider the spectral curve $p: C_a \to C$. 
    The companion Higgs bundle $(E_a = p_{\ast}\mc{O}_{C_a},\theta_a)$ obtains a fiberwise-nondegenerate symmetric pairing $\omega: S^2E_a \to (\Omega^1_C)^{-2n}$ that is anti-self-adjoint with respect to $\theta_a$. 
    By replacing $E_a$ with $E_a' = E_a \otimes (\Omega^1_C)^{n} = p_{\ast}p^{\ast}(\Omega^1_C)^{n}$, we obtain a $G$-Higgs bundle on $C$. 
    The determinant of $E_a'$ comes with a natural trivialization once a square root of $\Omega^1_C$ is chosen. 
\end{enumerate}
\begin{enumerate}[label=($C_n$)]
    \item \label{item:sympl_gp} $G = \rm{Sp}_{2n}$. 
    Recall that $A_{2n} \to A$ is given by the quotient of $(a_1,a_3,\ldots,a_{2n-1}) \subset A_{2n}$, so
    \[
        B = B_{2n} \otimes_{A_{2n}}A = A[x]/(x^{2n}+a_2x^{2n-2}+\cdots+a_{2n}), 
    \]
    Note that $B$ is a smooth $k$-algebra by the Jacobian criterion.

    There is a bilinear form $\omega: B \otimes_A B \to A$ that is nondegenerate, alternating, and anti-self-adjoint with respect to $\theta$. 
    Furthermore, $\omega$ satisfies the following $\G_m$-compatibility condition:
    \begin{equation*}
        \label{eq:G_m_compat_sympl_gp}
        (\omega \otimes \id)(\rm{coact}_{B\otimes_A B}(b_1\otimes b_2)) = \rm{coact}_A(\omega(b_1\otimes b_2))(1\otimes t^{2n-1}).
    \end{equation*}
    Let $a \in \mc{A}(C,G) \subset \mc{A}(C,\rm{GL}_{2n})$ and consider the spectral cover $p: C_a \to C$.
    The companion Higgs bundle $(E_a = p_{\ast}\mc{O}_{C_a},\theta_a)$ obtains a fiberwise, nondegenerate symplectic pairing $\omega: \Lambda^2 E_a \to (\Omega^1_C)^{1-2n}$ that is anti-self-adjoint with respect to $\theta$. 
    Choose a square root of $\Omega^1_C$ and denote it by $\Omega'$. 
    By replacing $E_a$ with $E_a' := E_a \otimes (\Omega')^{2n-1} = p_{\ast}p^{\ast}(\Omega')^{2n-1}$, we obtain a $G$-Higgs bundle on $C$ whose spectral data is $a$. 
\end{enumerate}
\begin{enumerate}[label=($D_n$)]
    \item \label{item:even_sp_orth_gp} $G = \rm{SO}_{2n}$. 
    Recall that $A = k[a_2,\ldots,a_{2n-2},p_n]$ with $p_n^2 = \det = a_{2n}$. 
    There is no longer a closed embedding $\mf{c} \to \mf{c}_{2n}$, and instead of working with $B = B_{2n} \otimes_{A_{2n}} A$, we work with 
    \[
        \tilde{B} = A[x,p_{n-1}]/(p_n-xp_{n-1},x^{2n-2}+a_2x^{2n-4}+\cdots + a_{2n-2} + p_{n-1}^2),
    \]
    which is the blowup of the singular algebra $B$ along the equation $(x)$.
    Observe that $\tilde{B}$ is a smooth $k$-algebra and a free $A$-module of rank $2n$. 

    Furthermore, $\tilde{B}$ is a flat $B$-module, as now explain.
    There is an injective $B$-module homomorphism $\tilde{B} \to B \otimes_A \tilde{B}$ that splits by \cite{atiyah_macdonald}*{Exercise 2.13}.
    The $B$-module $B \otimes_A \tilde{B}$ is free, because $\tilde{B}$ is a free $A$-module (of rank $2n$), so the $B$-module $\tilde{B}$ is projective, hence flat.
    
    Let $\theta: \tilde{B} \to \tilde{B}$ be the $A$-module endomorphism given by multiplication by $\ol{x}$. 
    Again, we can equip $\tilde{s} = \Spec(\tilde{B})$ with a $\G_m$-action that makes each map in the composition $\tilde{\mf{s}} \to \mf{s} \to \mf{c}$ equivariant. 

    There is a bilinear form $\omega: \tilde{B} \otimes_A \tilde{B} \to A$ that is nondegenerate, symmetric, and anti-self-adjoint with respect to $\theta$ with the following $\G_m$-compatibility condition:
    \begin{equation*}
        \label{eq:G_m_compat_even_spec_orth_gp}
        (\omega \otimes \id)(\rm{coact}_{B\otimes_A B}(b_1\otimes b_2)) = \rm{coact}_A(\omega(b_1\otimes b_2))(1\otimes t^{2n-2}).
    \end{equation*}

    Given $a \in \mc{A}(C,G)$, we obtain a finite, flat cover $\tilde{p}: \tilde{C}_a \to C_a$ via base change with $T^{\ast}_C(\tilde{\mf{s}}) \to T^{\ast}_C(\mf{c})$, which factors through $C_a$.
    Thus, the pushforward of $\mc{O}_{\tilde{C}_a}$ to $C$ is a $\rm{GL}_{2n}$-Higgs bundle with spectral data $a$. 
    Denote it by $(\tilde{E}_a,\tilde{\theta}_a)$ and observe that we have with a fiberwise nondegenerate symmetric pairing $\omega: S^2 \tilde{E}_a \to (\Omega^1_C)^{2-2n}$ that is anti-self-adjoint with respect to $\tilde{\theta}_a$. 
    By replacing $\tilde{E}_a$ with $\tilde{E}_a' := \tilde{E}_a \otimes (\Omega^1_C)^{n-1} = \tilde{p}_{\ast}\tilde{p}^{\ast}(\Omega^1_C)^{n-1}$, we obtain a $G$-Higgs bundle on $C$. 
    There is a natural trivialization of $\det(\tilde{E}_a')$ once we have chosen a square root of $\Omega^1_C$. 
\end{enumerate}

\subsubsection{Cohen-Macaulay spectral surfaces}\label{sect:CM_spectral_sf}
There is an analogous spectral correspondence for any smooth projective variety $X$ of dimension $d$ and $G = \rm{GL}_n$.
Under this correspondence, torsion-free Higgs sheaves $(E,\theta)$ on $X$ correspond to pure, $d$-dimensional coherent sheaves $\mc{E}$ on $T^{\ast}_X$ \cite{simpson_moduli2}*{Lemma 6.8}.
Moreover, a Higgs sheaf $(E,\theta)$ is $p$-semistable if and only if the corresponding sheaf $\mc{E}$ is a $p$-semistable coherent sheaf on $T^{\ast}_X$ \cite{simpson_moduli2}*{Corollary 6.9}.
Finally, there is a Hitchin-type morphism $\sigma$ from the moduli space of Higgs sheaves on $X$ to $\mc{A}(X,\rm{GL}_n)$ \cite{simpson_moduli2}*{\S Hitchin's proper map}, and the coherent sheaf $\mc{E}$ corresponding to a Higgs sheaf $(E,\theta)$ is scheme-theoretically supported on the (naive) spectral cover $X_{\sigma(E,\theta)} \subset T^{\ast}_X$ \cite{simpson_moduli2}*{Lemma 6.10 and following Remark}, \cite{chen_ngo20}*{Proposition 6.1}.

We are interested in \textit{locally-free} Higgs sheaves, i.e. Higgs bundles, and thus want to pushforward Cohen-Macaulay sheaves that are supported on a spectral cover to $X$.
However, a spectral cover may not be flat over $X$ as soon as $\dim(X) \geq 2$, so it is unclear how to find such Cohen-Macaulay sheaves.
To this end, for certain $b \in \mc{B}(X,\rm{GL}_n)$, Chen and Ng\^o construct a Cohen-Macaulayfication $X^{\rm{CM}}_b$ of the naive spectral cover $X_b$, and a finite, flat map $X^{\rm{CM}}_b \to X$. 
In particular, the pushforward of a line bundle on $X^{\rm{CM}}_b$ to $X$ is a Higgs bundle whose spectral data is $b$.
We elaborate on their construction in the subsequent paragraphs, and remark that the naive spectral cover is still useful when proving statements about the $p$-semistability of a Higgs bundle on $X$. 

In this subsection, we work with $G = \rm{GL}_n$ and a fixed smooth projective surface $X$, so we drop these symbols whenever it makes the notation neater. 
In this setting, the polarization ring $k[\mf{t}^2]^{\rm{pol}}$ is isomorphic to $k[\mf{t}^2]^W$ \cite{weyl_classical_groups}*{\S II.A.3}, and we can identify $\mf{B}$ with the Chow variety of $0$-cycles on $\A^2$ of degree $n$, denoted $\rm{Chow}_n(\A^2)$. This identifies $\mc{B}$ as the scheme of sections of $\rm{Chow}_n(T^{\ast}_X/X) \to X$. 
There is a $\rm{GL}_2$-invariant subscheme $Q \subset \rm{Chow}_n(\A^2)$, which is the complement of the locus of multiplicity-free $0$-cycles. 
Define the open subscheme $\mc{B}^{\heartsuit} \subset \mc{B}$ as the locus of points where $b: X \to \rm{Chow}_n(T^{\ast}_X/X)$ does not factor through $T^{\ast}_X(Q)$. 

For each $b \in \mc{B}^{\heartsuit}(k)$, Chen and Ng\^o construct a finite, flat covering $p_b^{\rm{CM}}: X_b^{\rm{CM}} \to X$ of degree $n$ such that the morphism $p_b^{\rm{CM}}$ factors through the naive spectral surface $X_b \subset T^{\ast}_X$, \cite{chen_ngo20}*{Proposition 7.2}. 
We call these covers CM-spectral surfaces. 

There is a spectral correspondence for these covers also proven by Chen and Ng\^o. 
Namely, for each $b \in \mc{B}^{\heartsuit}(k)$, the fiber of the Hitchin morphism $h^{-1}(b)$ is isomorphic to the stack of Cohen-Macaulay sheaves of generic rank 1 on $X_b^{\rm{CM}}$, \cite{chen_ngo20}*{Theorem 7.3}. 
Higgs bundles with spectral data $b$ are obtained by pushing forward Cohen-Macaulay sheaves of generic rank 1 on $X_b^{\rm{CM}}$. 
In particular, the Hitchin fiber $h^{-1}(b)$ contains the Picard stack of line bundles on the CM-spectral cover. 

\begin{lem} 
\label{lem:integral_CM_cover_implies_semistable}
    Let $b \in \mc{B}^{\heartsuit}(k)$. 
    Suppose that $X_b^{\rm{CM}}$ is integral, that $(E,\theta)$ is a Higgs bundle on $X$ with spectral data $b$, and that $E$ is the pushforward of a line bundle $\mc{E}$ on $X_b^{\rm{CM}}$. 
    Then $(E,\theta)$ is $p$-semistable.
    \begin{proof}
        We use the spectral correspondence between Higgs sheaves on $X$ and sheaves on $T^{\ast}_X$, which is proven in \cite{simpson_moduli2}*{Lemma 6.8}. 
        Under this correspondence, the $p$-semistability of a Higgs sheaf $(E,\theta)$ on $X$ is equivalent to the $p$-semistability of the corresponding coherent sheaf $\mc{E}$ on $T^{\ast}_X$ \cite{simpson_moduli2}*{Corollary 6.9}.
        
        Consider the composition
        \[
            \xymatrix{
                X_b^{\rm{CM}} \ar[r]^-{q} & X_b \ar[r]^-{p_b} & X.
            }
        \]
        If $F \subset E$ were a destabilizing sub-Higgs sheaf, then $F$ must equal $(p_b)_{\ast}\mc{F}$, where $\mc{F}$ is a destabilizing subsheaf of the rank one coherent sheaf $q_{\ast}\mc{E}$ on $X_b$. 
        Since $X_b^{\rm{CM}}$ is assumed to be integral, the sheaf $q_{\ast}\mc{E}$ is supported on $(X_b)_{\rm{red}}$, which is integral. 
        Thus there are no destabilizing subsheaves of $q_{\ast}\mc{E}$, and the lemma follows.
    \end{proof}
\end{lem}

\section{Pullback Higgs bundles}\label{sect:pullback_higgs_bdles}
\subsection{The case of a fibered surface}\label{sect:pullback_fibered_sf}
Let $f: X \to C$ be a fibered surface \ref{hyp:fs}. 
We show how the map of k\"ahler differentials $df: f^{\ast}\Omega^1_C \to \Omega^1_X$ induces morphisms between various spaces over $X$ and $C$ that appear in the definitions of the moduli of $G$-Higgs bundles, the Hitchin base, and the space of spectral data.

Suppose that $\{U_{\alpha} \to X\}_{\alpha \in I}$ is a finite cover by affine opens over which the locally free sheaves $f^{\ast}\Omega^1_C$ and $\Omega^1_X$ trivialize. 
Let $h_{\alpha\beta} \in \mc{O}(U_{\alpha\beta})^{\times}$ be the transition functions for $f^{\ast}\Omega^1_C$, and let $g_{\alpha\beta} \in \rm{GL}_2(\mc{O}(U_{\alpha\beta}))$ be the transition functions for $\Omega^1_X$. 
Then the cotangent morphism $f^{\ast}T^{\ast}_C \to T^{\ast}_X$ over $U_{\alpha}$ is given by an $R_{\alpha}$-algebra map
\begin{equation*}
    \begin{gathered}
        R_{\alpha}[z^1, z^2] \to R_{\alpha}[w] \\
        z^i \mapsto \phi^i_{\alpha}w,
    \end{gathered}
\end{equation*}
where $\phi^i_{\alpha} \in k$. Given any open $\Spec(R) \subset U_{\alpha\beta}$, we can extend scalars to $R$, and the collection of $\phi^i_{\alpha}$ must satisfy  
\begin{equation}
    \label{eq:cotangent_morphism_locally}
    h_{\alpha\beta}\phi^i_{\alpha} = (g_{\alpha\beta})^i_j \phi^j_{\beta} \in R.
\end{equation}

The spaces $f^{\ast}T^{\ast}_C(\mf{g})$ and $T^{\ast}_X(\mf{C}^2)$ are affine over $X$, and we construct a morphism between them by producing a morphism of the corresponding $\mc{O}_X$-algebras. 
Let $V$ be the $k$-vector space $k^2$. 
Identify $\mf{g}^2$ with $\Hom_k(V^{\ast},\mf{g})$. 
Let $e_1,e_2$ be a basis of $V$, and let $e^1,e^2$ be the dual basis. 
Let $f_1,\ldots,f_m$ be a basis of $\mf{g}$, let $f^1,\ldots,f^m$ be the dual basis, and let $\epsilon^{\ell}_{ij}$ denote the structure constants for $\mf{g}$ with respect to this basis. 

Under this identification, the $\rm{GL}_2$-action on $\mf{g}^2$ is given by the coaction map
\begin{gather*}
    \sigma: k[\mf{g}^2] \to k[\mf{g}^2]\otimes_k k[\rm{GL}_2] \\
    e^i \otimes f^{\ell} \mapsto (e^j\otimes f^{\ell})\otimes x^i_j.
\end{gather*}
This induces a coaction map on $k[\mf{C}^2]$, since this ring is the quotient of $k[\mf{g}^2]$ by the ideal $I$ generated by the elements $\epsilon^{\ell}_{ij}(e^1\otimes f^i)(e^2 \otimes f^j)$, and these generators are sent to zero by $\sigma$ composed with the map to $k[\mf{C}^2] \otimes k[\rm{GL}_2]$. 

Next, we take the ring map corresponding to the diagonal embedding $\mf{g} \to \mf{g}^2$, then tensor with $R_{\alpha}$ and twist by $df$: 
\begin{equation*}
\label{eq:diag_embedding}    
\begin{gathered}
    R_{\alpha}[e^1,e^2]\otimes_k k[\mf{g}] \to R_{\alpha} \otimes_k k[\mf{g}] \\
    (e^i \otimes f^j) \mapsto 1 \otimes \phi^i_{\alpha}f^j.
\end{gathered}
\end{equation*}
Observe that the generators of $I$ are sent to zero:
\begin{equation*}
\label{eq:ideal_killed}
    \epsilon^{\ell}_{ij}(e^1\otimes f^i)(e^2 \otimes f^j) \mapsto 1 \otimes \epsilon^{\ell}_{ij}(\phi^1_{\alpha}f^i)(\phi^2_{\alpha}f^j) 
    = 1 \otimes \phi^1_{\alpha}\phi^2_{\alpha}(\epsilon^{\ell}_{ij}f^if^j) 
    = 0,
\end{equation*}
since $\epsilon^{\ell}_{ij} = -\epsilon^{\ell}_{ji}$ for all $i,j,\ell$. Therefore we get induced maps $R_{\alpha} \otimes k[\mf{C}^2] \to R_{\alpha} \otimes k[\mf{g}]$. 

These morphisms are compatible with the transition maps. 
Given $\Spec(R) \subset U_{\alpha\beta}$, we have a commutative diagram
\begin{equation*}
\label{eq:compat_trans_maps}
    \xymatrix{
        R \otimes k[\mf{C}^2] \ar[r] \ar[d] & R\otimes k[\mf{g}] \ar[d] \\
        R \otimes k[\mf{C}^2] \ar[r] & R\otimes k[\mf{g}],
    }
    \hspace{1.5cm}
    \xymatrix{
        (e^i\otimes f^{\ell}) \ar@{|->}[r] \ar@{|->}[d] & 1 \otimes \phi^i_{\alpha}f^{\ell} \ar@{|->}[d] \\
        (g_{\alpha\beta})^i_je^j\otimes f^{\ell} \ar@{|->}[r] & (g_{\alpha\beta})^i_j \otimes \phi^j_{\beta}f^{\ell} = h_{\alpha\beta}\otimes \phi^i_{\alpha}f^{\ell}, 
    }
\end{equation*}
where the equality holds because of equation (\ref{eq:cotangent_morphism_locally}). 
Therefore we obtain an $X$-morphism from $f^{\ast}T^{\ast}_C(\mf{g})$ to $T^{\ast}_X(\mf{C}^2)$.

The group $G$ acts on $f^{\ast}T^{\ast}_C(\mf{g})$ and $T^{\ast}_X(\mf{C}^2)$ compatibly with their respective maps to $X$. 
The action is the fiberwise (simultaneous) adjoint action. 
The morphism we just constructed is $G$-equivariant, because each $R_{\alpha}\otimes k[\mf{C}^2] \to R_{\alpha} \otimes k[\mf{g}]$ is $G$-equivariant, and the $G$-action commutes with the action of $\rm{GL}_2$ on $\mf{C}^2$, respectively $\G_m$ on $\mf{g}$. 

Recall that in \S \ref{sect:hitchin_morphism}, we have a composition $\mf{C}^2 \to \mf{B} \to \mf{A}$, where each arrow is $\rm{GL}_2$-equivariant, and the map $\mf{C}^2 \to \mf{B}$ is $G$-invariant. 
Therefore for each $\alpha$, we have a commutative diagram of $R_{\alpha}$-algebras
\begin{equation*}
\label{eq:relating_affine_algs}
    \xymatrix{
        & R_{\alpha} \otimes k[\mf{A}] \ar[dr] \ar[dl] & \\
        R_{\alpha}\otimes k[\mf{B}] \ar[rr] \ar[d] & & R_{\alpha} \otimes k[\mf{C}^2] \ar[d] \\
        R_{\alpha} \otimes k[\mf{g}]^G \ar[rr] & & R_{\alpha} \otimes k[\mf{g}]. 
    }
\end{equation*}
Furthermore, these diagrams are compatible with the transition functions, so we have established the following
\begin{prop}
\label{prop:spaces_over_X}
Let $f: X \to C$ be a fibered surface \ref{hyp:fs}, and let $G$ be a reductive group.
\begin{enumerate}
    \item There is a commutative diagram of spaces over $X$:
        \begin{equation*}
            \label{diag:spaces_over_X}
            \xymatrix{
                [f^{\ast}T^{\ast}_C(\mf{g})/G] \ar[rr] \ar[d] & & [T^{\ast}_X(\mf{C}^2)/G] \ar[d] \\
                f^{\ast}T^{\ast}_C(\mf{c}) \ar[rr] \ar[dr] & & T^{\ast}_X(\mf{B}) \ar[dl] \\
                & T^{\ast}_X(\mf{A}). & 
            }
        \end{equation*}
    \item Taking the scheme of sections of the triangle in the diagram gives us a composition $\mc{A}(C,G) \to \mc{B}(X,G) \to \mc{A}(X,G)$, which is the closed embedding induced by pullback of symmetric differentials
        \[ H^0(C,S^{i}\Omega^1_C) \to H^0(X,S^i\Omega^1_X). \]
        In particular, the Hitchin base $\mc{A}(C,G)$ embeds into $\mc{B}(X,G)$.
\end{enumerate}
\end{prop}

There is a cartesian diagram
\begin{equation*}
    \xymatrix{
        [f^{\ast}T^{\ast}_C(\mf{g})/G] \ar[r] \ar[d] & [T^{\ast}_C(\mf{g})/G] \ar[d] \\
        X \ar[r]_{f} & C,
    }
\end{equation*}
(see \S \ref{sect:assoc_bdle}), so a $G$-Higgs bundle on $C$ given by a section $h_{E,\theta}: C \to [T^{\ast}_C(\mf{g})/G]$ yields a $G$-Higgs bundle on $X$ via the composition
\begin{equation*}
    \xymatrix{
        X \ar[r]^-{\simeq} & X \times_X C \ar[r] & [f^{\ast}T^{\ast}_C(\mf{g})/G] \ar[r] & [T^{\ast}_X(\mf{C}^2)/G].
    }
\end{equation*}
To be explicit, the pullback $G$-Higgs bundle on $X$ is the $G$-bundle $f^{\ast}E$ and Higgs field $f^{\ast}\theta$ given by the composition 
\begin{equation*}
    \xymatrixcolsep{1.5cm}\xymatrix{
        \mc{O}_X \ar[r]^-{f^{\ast}\theta} & \rm{ad}(f^{\ast}E) \otimes f^{\ast}\Omega^1_C \ar[r]^-{\id \otimes\, df} & \rm{ad}(f^{\ast}E) \otimes \Omega^1_X.
    }
\end{equation*}
This discussion coupled with Proposition \ref{prop:spaces_over_X} yields the following
\begin{cor}
    \label{cor:pullback_higgs_bdle_spec_data}
    Let $f: X \to C$ be a fibered surface \ref{hyp:fs}, and let $(E,\theta)$ be a $G$-Higgs bundle on $C$. Then the pullback $G$-Higgs bundle $(f^{\ast}E,f^{\ast}\theta)$ on $X$ satisfies $h_{X,G}(f^{\ast}E,f^{\ast}\theta) = h_{C,G}(E,\theta)$, after viewing $\mc{A}(C,G)$ as a closed subscheme of $\mc{B}(X,G)$. 
\end{cor}

\subsection{Change of group}\label{sect:change_of_gp}
Let $G$ be a classical group and let $G \to \rm{GL}_N$ be the standard linear representation. 
We can choose maximal tori for each group so that there is an induced map of Cartan subalgebras $\mf{t} \to \mf{t}_N$ that intertwines the $W$-action on $\mf{t}$ with the $\mf{S}_N$-action on $\mf{t}_N$. 
Pick generators $a_1,\ldots,a_N$ of $k[\mf{t}_N]^{\mf{S}_N}$ and generators $c_1,\ldots,c_n$ of $k[\mf{t}]^W$. 
The restriction of $a_j$ to $\mf{t} \subset \mf{t}_N$ is some polynomial in the variables $c_1,\ldots,c_n$, which we denote by $f_j$.

Recall that the polarization ring $k[(\mf{t}_N)^2]^{\rm{pol}} \subset k[(\mf{t}_N)^2]^{\mf{S}_N}$ is generated by the $\ul{i}$-polarizations of $a_j$, for $j = 1,\ldots,N$ and $\ul{i}$ a weak composition of $j$ of length 2. 
The $\ul{i}$-polarization of $a_j$ restricted to $\mf{t}^2 \subset (\mf{t}_N)^2$ equals the $\ul{i}$-polarization of $f_j$, which defines a ring map $k[\mf{t}^2]^{\rm{pol}} \to k[(\mf{t}_N)^2]^{\rm{pol}}$. 
Furthermore, we can restrict polarizations to the image of the diagonal embedding $\mf{t}_N \to (\mf{t}_N)^2$, which gives us $\mf{S}_N$-invariant functions on $\mf{t}_N$. 
Restricting a polarization first to $\mf{t}_N$, then to $\mf{t}$ is the same as first restricting to $\mf{t}^2$, then to $\mf{t}$. 

Next, recall that the ring $k[\mf{A}_N]$ is a polynomial ring generated by symbols $z_{j,\ul{i}}$, and the ring map $k[\mf{A}_N] \to k[(\mf{t}_N)^2]^{\rm{pol}}$ is defined by sending $z_{j,\ul{i}}$ to $a_{j,\ul{i}}$ (respectively for $G$, send a generator $w_{j,\ul{i}} \in k[\mf{A}]$ to $c_{j,\ul{i}}$). 
Since the polarizations $f_{j,\ul{i}}$ are polynomials in terms of the $c_{j,\ul{i}}$, we can define a map $k[\mf{A}_N] \to k[\mf{A}]$, making the following diagram commutative:
\begin{equation*}
    \label{eq:chg_gp_diagram}
    \xymatrixcolsep{1.5cm}\xymatrix{
        k[\mf{A}_N] \ar[r] \ar[d] & k[(\mf{t}_N)^2]^{\rm{pol}} \ar[r]^-{\rm{restr}} \ar[d] & k[\mf{t}_N]^{\mf{S}_N} \ar[d] \ar@{^{(}->}[r] & k[\mf{t}_N] \ar[d] \\
        k[\mf{A}] \ar[r] & k[\mf{t}^2]^{\rm{pol}} \ar[r]_-{\rm{restr}} & k[\mf{t}]^W \ar@{^{(}->}[r] & k[\mf{t}].
    }
\end{equation*}

The left square is $\rm{GL}_2$-equivariant. 
The argument that produces a morphism $f^{\ast}T^{\ast}_C(\mf{g}) \to T^{\ast}_X(\mf{C}^2)$ shows (mutatis mutandis) that there is a commutative diagram of schemes over $X$: 
\begin{equation*}
\label{eq:chg_of_gp_schemes_over_X}
    \xymatrix{
        f^{\ast}T^{\ast}_X(\mf{c}) \ar[r] \ar[d] & T^{\ast}_X(\mf{B}) \ar[r] \ar[d] & T^{\ast}_X(\mf{A}) \ar[d] \\
        f^{\ast}T^{\ast}_X(\mf{c}_N) \ar[r] & T^{\ast}_X(\mf{B}_N) \ar[r] & T^{\ast}_X(\mf{A}_N).
    }
\end{equation*}
By taking scheme of sections, we obtain the following
\begin{prop}
\label{prop:change_of_group}
    Let $G$ be a classical group and let $G \to \rm{GL}_N$ be the standard representation. Then there is a commutative diagram relating the Hitchin bases and spaces of spectral data for $G$ and $\rm{GL}_N$:
    \begin{equation}
    \label{eq:relating_hitch_base_and_spec_data}
        \xymatrix{
            \mc{A}(C,G) \ar[r] \ar[d] & \mc{B}(X,G) \ar[r] \ar[d] & \mc{A}(X,G) \ar[d] \\
            \mc{A}(C,\rm{GL}_N) \ar[r]_{\iota_f} & \mc{B}(X,\rm{GL}_N) \ar[r] & \mc{A}(X,\rm{GL}_N).
        }
    \end{equation}
\end{prop}

\subsection{The case of a morphism between two surfaces}
Let $\pi: Y \to X$ be a morphism between two smooth, projective surfaces. 
Once again, we show how the map of k\"ahler differentials $d\pi: \pi^{\ast}\Omega^1_X \to \Omega^1_Y$ induces morphisms between various spaces over $Y$ and $X$ that appear in the definitions of the moduli of $G$-Higgs bundles, the Hitchin base, and the space of spectral data. 

Let $\{U_{\alpha} \to Y\}_{\alpha \in I}$ be an open affine cover over which both locally free sheaves $\pi^{\ast}\Omega^1_X$ and $\Omega^1_Y$ trivialize. 
Let $g_{\alpha\beta},g'_{\alpha\beta} \in \rm{GL}_2(\mc{O}(U_{\alpha\beta})$ denote the transition functions for $\Omega^1_Y$ and $\pi^{\ast}\Omega^1_X$ respectively.
Then the cotangent morphism $\pi^{\ast}T^{\ast}_X \to T^{\ast}_Y$ over $U_{\alpha}$ is given by an $R_{\alpha}$-algebra map
\begin{equation*}
    \begin{gathered}
        R_{\alpha}[z^1,z^2] \to R_{\alpha}[z^1,z^2] \\
        z^i \mapsto (\phi_{\alpha})^i_j z^j,
    \end{gathered}
\end{equation*}
where $(\phi_{\alpha})^i_j \in k$. 
Given any $V = \Spec(R) \subset U_{\alpha\beta}$, the collection of morphisms must satisfy
\begin{equation}
    \label{eq:compat_w_trans_fns_two_sfs}
    (g_{\alpha\beta})^i_j (\phi_{\beta})^j_k = (\phi_{\alpha})^i_j(g'_{\alpha\beta})^j_k \in R,
\end{equation}
for all $i,j,\alpha,\beta$. 

We use the same notation from \S \ref{sect:pullback_fibered_sf} in the constructions that follow.
For each $\alpha \in I$, define a map of $R_{\alpha}$-algebras
\begin{equation*}
    \begin{gathered}
        \Phi_{\alpha}: R_{\alpha}[e^1,e^2] \otimes_k k[\mf{g}] \to R_{\alpha}[e^1,e^2] \otimes_k k[\mf{g}] \\
        (e^i\otimes f^{\ell}) \mapsto (\phi_{\alpha})^i_j e^j \otimes f^{\ell}.
    \end{gathered}
\end{equation*}
The generators of the ideal $I \subset k[\mf{g}^2]$ corresponding to $\mf{C}^2$ are sent to zero by the map above.
Furthermore, these morphisms are compatible with the transition functions:
\begin{equation*}
    \xymatrixcolsep{2cm}\xymatrix{
        e^i \otimes f^{\ell} \ar@{|->}[r] \ar@{|->}[d] & (\phi_{\alpha})^i_j e^j \otimes f^{\ell} \ar@{|->}[d] \\
        (g_{\alpha\beta})^i_j e^j\otimes f^{\ell} \ar@{|->}[r] & (g_{\alpha\beta})^i_j(\phi_{\beta})^j_p e^p \otimes f^{\ell} = (\phi_{\alpha})^i_j (g'_{\alpha\beta})^j_p e^p \otimes f^{\ell},
    }
\end{equation*}
because of equation (\ref{eq:compat_w_trans_fns_two_sfs}). 
Once again, the group $G$ acts fiberwise on the two $Y$-schemes $\pi^{\ast}T^{\ast}_X(\mf{C}^2)$ and $T^{\ast}_Y(\mf{C}^2)$, and the morphism we have just constructed between them is $G$-equivariant. 

At this point, for each $\alpha \in I$, we have a commutative diagram of $R_{\alpha}$-algebras:
\begin{equation*}
    \xymatrix{
        R_{\alpha} \otimes_k k[\mf{A}] \ar@{-->}[d] \ar[r] & R_{\alpha} \otimes_k k[\mf{B}] \ar@{^{(}->}[r] \ar@{-->}[d] & R_{\alpha} \otimes_k k[\mf{t}^2]^W \ar[r]^{\simeq} \ar[d] & R_{\alpha} \otimes_k k[\mf{C}^2]^G \ar@{^{(}->}[r] \ar[d] & R_{\alpha} \otimes_k k[\mf{C}^2] \ar[d]^{\Phi_{\alpha}} \\
        R_{\alpha} \otimes_k k[\mf{A}] \ar[r] & R_{\alpha} \otimes_k k[\mf{B}] \ar@{^{(}->}[r] & R_{\alpha} \otimes_k k[\mf{t}^2]^W \ar[r]^{\simeq} & R_{\alpha} \otimes_k k[\mf{C}^2]^G \ar@{^{(}->}[r] & R_{\alpha} \otimes_k k[\mf{C}^2] 
    }
\end{equation*}
We can fill in the second dotted arrow once we verify that a polarization $p \in k[\mf{t}^2]^{\rm{pol}}$ is sent to $R_{\alpha} \otimes_k k[\mf{t}^2]^{\rm{pol}} \subset R_{\alpha} \otimes_k k[\mf{t}^2]^W$. 
To this end, let $(b_1,b_2) \in \A^2(k),(\theta^1,\theta^2) \in \mf{t}^2(k),(x^i_j) \in \rm{Mat}_2(k)$, and let $c_m \in k[\mf{t}]^W$ be a generator. 
Then we may compute $c_m(b_ix^i_j\theta^j)$ in two ways:
\begin{equation*}
    c_m(b_ix^i_j\theta^j) = \sum_{\ul{i}} b^{\ul{i}} c_{m,\ul{i}}(x^1_j\theta^j,x^2_j\theta^j) = \sum_{\ul{\ell}} (b_ix^i_1)^{\ell_1}(b_ix^i_2)^{\ell_2}c_{m,\ul{\ell}}(\theta^1,\theta^2).
\end{equation*}
By comparing the coefficient of $b^{\ul{i}}$ in each expression, we verify the claim and thus can fill in the second dotted arrow.

Recall that $k[\mf{A}]$ is a polynomial ring generated by the symbols $z_{j,\ul{i}}$ where $j$ ranges over $\{1,\ldots,r\}$ and $\ul{i}$ ranges over weak compositions of $e_j$ of length 2. 
The $R_{\alpha}$-algebra maps $R_{\alpha} \otimes k[\mf{A}] \to R_{\alpha} \otimes k[\mf{B}]$ in both rows are given by sending $z_{j,\ul{i}}$ to $c_{j,\ul{i}}$, and thus we can fill in the first dotted arrow making the entire diagram commutative. 

The composition $\mf{C}^2 \to \mf{B} \to \mf{A}$ (\S \ref{sect:hitchin_morphism}) is $\rm{GL}_2$-equivariant, so the above diagram is compatible with transition maps. 
Therefore we have established the following 
\begin{prop}
    \label{prop:diagram_morphism_of_sfs}
    Let $\pi: Y \to X$ be a morphism between two smooth, projective surfaces and let $G$ be a reductive group.
    Then there is a commutative diagram of spaces over $Y$
    \[
        \xymatrix{
            [\pi^{\ast}T^{\ast}_X(\mf{C}^2)/G] \ar[r] \ar[d] & \pi^{\ast}T^{\ast}_X(\mf{B}) \ar[r] \ar[d] & \pi^{\ast}T^{\ast}_X(\mf{A}) \ar[d] \\
            [T^{\ast}_Y(\mf{C}^2)/G] \ar[r] & T^{\ast}_Y(\mf{B}) \ar[r] & T^{\ast}_Y(\mf{A}).
        }
    \]
    By taking the stack of sections, we obtain a commutative diagram
    \[
        \xymatrix{
            \mc{M}(X,G) \ar[r]^{h_{X,G}} \ar[d]_{\pi^{\ast}} & \mc{B}(X,G) \ar[r] \ar[d] & \mc{A}(X,G) \ar[d] \\
            \mc{M}(Y,G) \ar[r]_{h_{Y,G}} \ar[r] & \mc{B}(Y,G) \ar[r] & \mc{A}(Y,G).
        }
    \]
\end{prop}

\subsection{Pullback of symmetric differentials for blowups}\label{sect:pullback_symm_diff_sfs}
Let $X$ be a smooth, projective surface, and let $\pi: Y \to X$ be the blowup of a point in $p \in X$. 
Let $j: U \to X$ denote the open immersion of the complement of $p$ in $X$, and let $j': U' \simeq U \to Y$ be the base change. 
By Proposition \ref{prop:diagram_morphism_of_sfs}, the pullback of symmetric differentials $H^0(X,S^i\Omega^1_X) \to H^0(Y,S^i\Omega^1_Y)$, which is injective for each $i$, induces a commutative diagram
\begin{equation*}
    \xymatrix{
        \mc{B}(X,G) \ar@{^{(}->}[r] \ar@{^{(}->}[d] & \mc{B}(Y,G) \ar@{^{(}->}[d] \\
        \mc{A}(X,G) \ar@{^{(}->}[r] & \mc{A}(Y,G),
    }
\end{equation*}
where each arrow is a closed immersion.

For each $i$, the pullback map $H^0(X,S^i\Omega^1_X) \to H^0(Y,S^i\Omega^1_Y)$ is injective.
A result of Serre implies that the adjoint pair $(j^{\ast},j_{\ast})$ induces an equivalence of categories between locally free sheaves on $X$ and locally free sheaves on $U$ \cite{serre_prolongement_de_faisceaux}*{Proposition 7}.
In particular, we have natural isomorphisms $S^i\Omega^1_X \isomto j_{\ast}S^i\Omega^1_U$ for each $i \geq 0$. 
Therefore, we get an injection
\[
    H^0(Y,S^i\Omega^1_Y) \hookrightarrow H^0(U',S^i\Omega^1_{U'}) \isomto H^0(U,S^i\Omega^1_U) \isomto H^0(X,S^i\Omega^1_X),
\]
from which it follows that $\mc{A}(X,G) = \mc{A}(Y,G)$.

\section{The Hitchin morphism for fibered surfaces and classical groups}\label{sect:hitch_morphism_fibered_sf_classical_gp}
\subsection{Relationship between different notions of spectral cover}\label{sect:relating_spec_covers}
Let $f: X \to C$ be a fibered surface \ref{hyp:fs}, let $G$ be a classical group, and let $G \to \rm{GL}_N$ be the standard linear representation. 
Given $a \in \mc{A}(C,\rm{GL}_N) \cap \mc{B}(X,\rm{GL}_N)^{\heartsuit}$, we can form a finite, flat cover $X \times_C C_a \to X$ (\S \ref{sect:spectral_covers}), and ask how it relates to $X_a^{\rm{CM}}$ (\S \ref{sect:CM_spectral_sf}). 
If we suppose that $f: X \to C$ has only reduced fibers, then for every such $a$, the cover $X \times_C C_a$ is isomorphic to $X_a^{\rm{CM}}$ \cite{chen_ngo20}*{Lemma 8.1}. 

\begin{lem}
\label{lem:integral_spectral_curve_implies_irreducible_base_change_X}
    If $C_a$ is integral, then $X \times_C C_a$ is irreducible.
    \begin{proof}
        The morphism $X \times_C C_a \to C_a$ is proper and flat by base change.
        It suffices to show that the generic fiber of this morphism is irreducible \cite{ega_4-2}*{Corollaire 2.3.5 (iii)}.
        To this end, consider the restriction of this morphism to the non-empty open subscheme $U \subset C_a$ where the fibers are geometrically reduced. 
        The function that sends a point $u \in U$ to the number of irreducible components of the geometric fiber over $u$ is upper semicontinuous \cite{landesman2016interpolationalgebraicgeometry}*{Proposition 3.2.5}.
        There exists some $u \in U$ such that the geometric fiber is irreducible, (simply take $u$ to be a point where the covering $C_a \to C$ is \'etale), so the locus in $C_a$ where the geometric fibers of the morphism $X \times_C C_a \to C_a$ are irreducible is open and non-empty, hence dense in $C_a$.
    \end{proof}
\end{lem}
Furthermore, if $C_a$ is smooth, then the CM-spectral cover $X_a^{\rm{CM}}$ is normal, hence integral \cite{chen_ngo20}*{Corollary 8.3}. 

We record two ancillary lemmas before proving that certain $G$-Higgs bundles on $X$ that are pulled back from $C$ are semiharmonic.
Recall that the open locus $\mc{A}(C,\rm{GL}_N)^{\rm{grss}} \subset \mc{A}(C,\rm{GL}_N)$, (where ``grss'' stands for ``generically regular semisimple''), is the locus where the spectral cover $C_a \to C$ is generically \'etale, and that we have a chain of open embeddings
\begin{equation*}
    \mc{A}(C,\rm{GL}_N)^{\rm{sm}} \subset \mc{A}(C,\rm{GL}_N)^{\rm{int}} \subset \mc{A}(C,\rm{GL}_N)^{\rm{grss}}.
\end{equation*}
\begin{lem}
\label{lem:grss_maps_to_heart}
    Let $G$ be a classical group and let $G \to \rm{GL}_N$ be the standard linear representation. 
    If $a \in \mc{A}(C,G)$ maps to $\mc{A}(C,\rm{GL}_N)^{\rm{grss}}$ in the diagram (\ref{eq:relating_hitch_base_and_spec_data}), then $a$ maps to $\mc{B}(X,\rm{GL}_N)^{\heartsuit} \subset \mc{B}(X,\rm{GL}_N)$. 
    \begin{proof}
        There exists an open, dense set $U_a \subset C$ such that for all $c \in U_a$, $a(c) = [\lambda_1,\ldots,\lambda_N] \in \rm{Chow}_N(T^{\ast}_{C,c})$ consists of $N$ distinct points. 
        Let $U \subset C$ be the largest open set over which $f: X \to C$ is smooth. 
        Then for every $x \in X \times_C (U\cap U_a)$, the 0-cycle $\iota_f(a)(x)$ consists of $N$ distinct points of $T^{\ast}_{X,x}$, where $\iota_f$ is the map defined in diagram (\ref{eq:relating_hitch_base_and_spec_data}).
    \end{proof}
\end{lem}
Lemma \ref{lem:grss_maps_to_heart} is used implicitly throughout the subsequent subsections.
The lemma essentially says that it makes sense to (define and) compare the covers $X \times_C C_a$ and $X_a^{\rm{CM}}$ whenever $a$ lies in $\mc{A}(C,\rm{GL}_N)^{\rm{grss}}$.

\begin{lem}
\label{lem:glob_gen_ev_is_smooth}
    Let $X$ be a smooth, connected, projective variety, and let $E$ be a locally free sheaf of rank $r$ on $X$. 
    Suppose that $E$ is globally generated. 
    Then the evaluation morphism
    \[
        \rm{ev}: X \times H^0(X,E) \to \uSpec(S^{\bullet}E^{\vee}),
    \]
    is smooth.
    \begin{proof}
        Observe that both the source and target are smooth over $k$, so it suffices to show for every closed point $(x,s) \in X \times H^0(X,E)$, the induced map on Zariski tangent spaces is surjective \cite{hart77}*{Proposition III.10.4}.

        The tangent space of $X \times H^0(X,E)$ at $(x,s)$ is $T_{X,x} \times H^0(X,E)$, and the tangent space of $\uSpec(S^{\bullet}E^{\vee})$ at $\rm{ev}(x,s)$ is $T_{X,x} \times E(x)$. 
        The map $d(\rm{ev})_{x,s}$ is given by the identity map on the $T_{X,x}$-factor, and the natural evaluation map $H^0(X,E) \to E(x)$ on the second factor. 
        This map is surjective for all $x$ because $E$ is globally generated. 
    \end{proof}
\end{lem}

We use the notations from \S\ref{sect:companion_sect_classical_gp} and assume that $g(C) \geq 1$ for each subsequent subsection.
For each classical group $G$, and for generic $a \in \mc{A}(C,G)$, we show that the pullback of the companion $G$-Higgs bundle $(E_a',\theta_a')$ is a semiharmonic $G$-Higgs bundle on $X$. 

\subsection{The special linear group}\label{sect:spec_lin_gp}
Assume that $g(C) \geq 2$, and apply Lemma \ref{lem:glob_gen_ev_is_smooth} to $E = \oplus_{i = 2}^{n+1} (\Omega^1_C)^{i}$ to see that the evaluation map $C \times \mc{A}(C,G) \to T^{\ast}_C(\mf{c})$ is smooth. 
If $C$ is a smooth elliptic curve, then $T^{\ast}_C$ is a trivial bundle of rank 1, and thus the evaluation morphism $C \times \mc{A}(C,G) \to T^{\ast}_C(\mf{c})$ may be identified with the identity morphism of $C \times \A^n$. 
In particular, the evaluation morphism is smooth in this case.

Because $T^{\ast}_C(\mf{s})$ is smooth over $\Spec(k)$, so too is the relative curve
\[ C_{\mc{A}} = (C \times \mc{A}(C,G)) \times_{T^{\ast}_C(\mf{c})} T^{\ast}_C(\mf{s}).\]
Apply generic smoothness to the map $C_{\mc{A}} \to \mc{A}(C,G)$ to see that the spectral curve $C_a$ is smooth for generic $a \in \mc{A}(C,G)$. 

Thus for such $a \in \mc{A}(C,G) \subset \mc{A}(C,\rm{GL}_{n+1})$, the CM-spectral cover $X \times_C C_a \cong X_a^{\rm{CM}}$ is normal. 
The companion $G$-Higgs bundle on $C$ pulled back to $X$ is isomorphic to the pushforward of a line bundle from $X_a^{\rm{CM}}$ by flat base change. 
Since $X_a^{\rm{CM}}$ is integral, this Higgs bundle is $p$-semistable by Lemma \ref{lem:integral_CM_cover_implies_semistable}, and the Chern classes vanish because the determinant line bundle is trivial and because the bundle is pulled back from a curve. 

\subsection{The symplectic group}\label{sect:sympl_gp}
The argument from \S \ref{sect:spec_lin_gp} for the special linear group applies mutatis mutandis to the symplectic group. 

\subsection{The even special orthogonal group}\label{sect:even_sp_orth_gp}
Assume that $g(C) \geq 2$. Apply Lemma \ref{lem:glob_gen_ev_is_smooth} to $E = \oplus_{i = 1}^{n} (\Omega^1_C)^{2i} $ to see that the evaluation map is smooth. 
If $C$ is a smooth elliptic curve, then $T^{\ast}_C$ is a trivial bundle of rank 1, and thus the evaluation morphism $C \times \mc{A}(C,G) \to T^{\ast}_C(\mf{c})$ may be identified with the identity morphism of $C \times \A^n$. 
In particular, the evaluation morphism is smooth in this case.

Because $T^{\ast}_C(\tilde{s})$ is smooth over $k$, so too is the relative curve 
\[ \tilde{C}_{\mc{A}} = (C \times \mc{A}(C,G)) \times_{T^{\ast}_C} T^{\ast}_C(\tilde{s}). \]
Apply generic smoothness to $\tilde{C}_{\mc{A}} \to \mc{A}(C,G)$, to see that the curve $\tilde{C}_a$ is smooth for generic $a \in \mc{A}(C,G)$. 

Recall that $B \to \tilde{B}$ is flat. 
If $\tilde{C}_a$ is smooth, then faithfully flat descent implies that $C_a$ is integral, so $X \times_C C_a$ is isomorphic to $X_a^{\rm{CM}}$. 
The argument of \cite{chen_ngo20}*{Corollary 8.3} shows that $X \times_C \tilde{C}_a$ is normal, so $X_a^{\rm{CM}}$ is integral. 
We conclude that the pullback of the companion $G$-Higgs bundle from $C$ to $X$ is semiharmonic by adapting the last part of the argument in \S \ref{sect:spec_lin_gp}.

\subsection{The odd special orthogonal group}\label{sect:odd_sp_orth_gp}
For any $a \in \mc{A}(C,G)$, the spectral cover $C_a$ has a component isomorphic to the zero section of $C$ in $T^{\ast}_C$. 
One can see this either from the definition $B := A[x]/(x(x^{2n}+a_2x^{2n-2} + \cdots + a_{2n}))$, or from the fact that the companion $G$-Higgs bundle $(E_a',\theta_a')$ has a non-trivial sub-Higgs bundle isomorphic to $(\Omega^1_C)^{-n}$ contained in the kernel of $\theta_a'$. 

Let $p: T^{\ast}_C \to C$ be the projection, let $t$ be the tautological section of $p^{\ast}\Omega^1_C$, and let $a = (a_2,\ldots,a_{2n}) \in \oplus_{i=1}^n H^0(C,(\Omega^1_C)^{2i})$.
The curve $C_a$ is the zero scheme of the section of $p^{\ast}(\Omega^1_C)^{2n+1}$ given by $t(t^{2n} + p^{\ast}a_2t^{2n} + \cdots + p^{\ast}a_{2n})$.
The argument of \S \ref{sect:spec_lin_gp} implies that for generic $a$, the spectral curve $C_a$ decomposes into two irreducible components $C_a = Z \cup C_a'$, where $Z$ is the zero-section of $T^{\ast}_C$, and $C_a'$ is a smooth irreducible curve. 

\begin{lem}
\label{lem:C_a_prime_smooth_implies}
    If $C_a'$ is smooth, then $a \in \mc{A}(C,G) \cap \mc{A}(C,\rm{GL}_{2n+1})^{\rm{grss}}$. 
    \begin{proof}
        If we assume that $C_a'$ is smooth, then there exists a nonempty open subset $U \subset C$ such that $a(x) \in T^{\ast}_{C,x}$ consists of $2n$ distinct points for all $x \in U$. 
        Furthermore, we must have $a_{2n} \neq 0 \in H^0(C,(\Omega^1_C)^{2n})$, so there is a nonempty open set $U'\subset U$ such that $a(x) \in T^{\ast}_{C,x}$ does not intersect $0$ for all $x \in U'$. 
        Then over $U'$, the map $C_a \to C$ is unramified, and the lemma follows.
    \end{proof}
\end{lem}
\begin{lem}
\label{lem:two_pieces}
    Suppose $a \in \mc{A}(C,G)$ is chosen so that the spectral curve $C_a$ decomposes into two irreducible components $C_a = Z \cup C_a'$, and that $C_a'$ is smooth.
    Then the CM-spectral cover $X_a^{\rm{CM}}$ also breaks into two pieces, namely $X \times_C Z$, which is isomorphic to $X$, and $X \times_C C_a'$.
    \begin{proof}
        By Lemma \ref{lem:C_a_prime_smooth_implies}, we have $X \times_C C_a \cong X_a^{\rm{CM}}$.
        Every irreducible component $T \subset X\times_C C_a$ surjects onto an irreducible component of $C_a$ by \cite{ega_4-2}*{Corollaire 2.3.5 (ii)} and the fact that $X \times_C C_a \to C_a$ is proper.
        The proof of Lemma \ref{lem:integral_spectral_curve_implies_irreducible_base_change_X} implies that $X \times_C Z$ and $X \times_C C_a'$ are irreducible components of $X\times_C C_a$, so these are all the irreducible components of $X \times_C C_a \cong X_a^{\rm{CM}}$.
    \end{proof}
\end{lem}

The only obstruction for the pullback companion $G$-Higgs bundle $f^{\ast}E_a'$ to be semiharmonic is $p$-semistability, because the Chern classes vanish in rational cohomology. 
This vanishing implies that $p$-semistability for $f^{\ast}E_a'$ is equivalent to $\mu$-semistability \cite{simpson_moduli2}*{Remark following Corollary 6.7} and \cite{huybrechts_lehn}*{Corollary 1.6.9}. 
As is true for vector bundles on curves, $f^{\ast}E_a'$ is $\mu$-semistable if and only if its tensor product with any line bundle on $X$ is $\mu$-semistable. 
In particular, it suffices to check the $\mu$-semistability of $f^{\ast}E_a' \otimes f^{\ast}(\Omega^1_C)^{-n} \cong f^{\ast}E_a$, which is the pushforward of $\mc{O}$ from $X_a^{\rm{CM}}$ to $X$ by flat base change.

Once again, we use the spectral correspondence between Higgs sheaves on $X$ and coherent sheaves on $T^{\ast}_X$, which preserves $\mu$-semistability \cite{simpson_moduli2}*{Lemma 6.8, Corollary 6.9}.
It suffices to show that no \textit{saturated} sub-Higgs sheaves of $f^{\ast}E_a$ destabilize it \cite{simpson_moduli1}*{p. 89}.
The quotient of $f^{\ast}E_a$ by a proper saturated subsheaf is a proper quotient sheaf that is pure, hence torsion-free. 
Therefore by \cite{stacks-project}*{\href{https://stacks.math.columbia.edu/tag/0EBG}{Tag 0EBG},\href{https://stacks.math.columbia.edu/tag/0B3N}{Tag 0B3N}}, it suffices to show that there are no locally free sub-Higgs sheaves of $f^{\ast}E_a$ that destabilize it. 

This discussion shows that the $\mu$-semistability of $f^{\ast}E_a$ as a Higgs bundle is equivalent to the $\mu$-semistability of $\mc{O}_{X_a^{\rm{CM}}}$ as a coherent sheaf by the CM-spectral correspondence \cite{chen_ngo20}*{Theorem 7.3}.
The only quotient sheaves that can destabilize $\mc{O}_{X_a^{\rm{CM}}}$ are the pushforwards to $X_a^{\rm{CM}}$ of $\mc{O}_{X\times_C Z}$ and $\mc{O}_{X\times_C C_a'}$, as these are the only proper quotient sheaves of pure dimension 2. 
Because of the spectral correspondence and the fact that $X_a^{\rm{CM}} \to X$ is finite, we can check the relevant inequalities on $X$.
Note that $\mc{O}_{X \times_C Z}$ pushes forward to $\mc{O}_X$, and $\mc{O}_{X\times_C C'_a}$ pushes forward to $\oplus_{i = 0}^{2n-1}f^{\ast}(\Omega^1_C)^{-i}$.

Fix an ample line bundle $L$ on $X$, and observe that the following inequalities hold
\begin{equation*}
    \begin{aligned}
        \mu(f^{\ast}E_a) = \frac{c_1(f^{\ast}E_a)\cdot L}{2n+1} = \frac{-(2n)(2n+1)c_1(f^{\ast}\Omega^1_C)\cdot L}{2n+1} &= (-n)(c_1(f^{\ast}K_C)\cdot L) \\
        &\leq 0,
    \end{aligned}
\end{equation*}
and 
\begin{equation*}
    \begin{aligned}
        \mu(f^{\ast}E_a) = \frac{c_1(f^{\ast}E_a)\cdot L}{2n+1} = \frac{-(2n)(2n+1)c_1(f^{\ast}\Omega^1_C)\cdot L}{2n+1} &= (-n)(c_1(f^{\ast}\Omega^1_C)\cdot L) \\
        &\leq \frac{-(2n-1)}{2}(c_1(f^{\ast}\Omega^1_C)\cdot L),
    \end{aligned}
\end{equation*}
because the quantity $c_1(f^{\ast}\Omega_C)\cdot L$ is positive. 

\phantomsection
\addcontentsline{toc}{section}{References}
\bibliography{references}

{\footnotesize
    \textsc{Department of Mathematics, Stony Brook University,
        Stony Brook, NY 11794-3651,
    USA}\par\nopagebreak
      \textit{E-mail address}: \texttt{matthew.huynh@stonybrook.edu}
}
\end{document}